\documentclass[12pt,psamsfonts]{amsart}
\usepackage{amssymb,amsmath,amscd,eucal,epsfig}
\usepackage{tikz}

 \def\Dj{\hbox{D\kern-.73em\raise.30ex\hbox{-} \raise-.30ex\hbox{}}}
 \def\dj{\hbox{d\kern-.33em\raise.80ex\hbox{-} \raise-.80ex\hbox{\kern-.40em}}}

\def\<{\langle}                     
\def\>{\rangle}                     

\newtheorem{thm}{Theorem}[section]

\newtheorem{lem}[thm]{Lemma}

\newcommand{\ben}{\begin{enumerate}}
\newcommand{\een}{\end{enumerate}}

\theoremstyle{plain}

\newtheorem{corollary}{Corollary}[section]

\theoremstyle{definition}

\newtheorem{rem}{Remark}[section]

\numberwithin{equation}{section}
\begin{document}
\title[Maximmum Induced Matching Numbers of Grids]{Some bounds on the maximum induced matching numbers of certain grids}
\author[D.O. Ajayi]{Deborah Olayide Ajayi$^1$}
\address{$^1$Department of Mathematics,
\newline \indent University of Ibadan,
\newline \indent Ibadan,
\newline \indent Nigeria}
\email{olayide.ajayi@mail.ui.edu.ng; adelaideajayi@yahoo.com}

\author[T.C Adefokun]{Tayo Charles Adefokun$^2$ }
\address{$^2$Department of Computer and Mathematical Sciences,
\newline \indent Crawford University,
\newline \indent Nigeria}
\email{tayoadefokun@crawforduniversity.edu.ng; tayo.adefokun@gmail.com}


\keywords{Induced Matching, Grids, Maximum Induced Matching Number
\\
\indent 2010 {\it Mathematics Subject Classification}. Primary: 05C70, 05C15}
%
%
%
\begin{abstract}
An induced matching $M$ in a graph $G$ is a matching in $G$ that is also the edge set of an induced subgraph of $G$. That is, any edge not in $M$ must have no more than one incident vertex saturated by $M$. The maximum size $|M|$ of an induced matching $M$ of $G$ is maximum induced matching number of $G$, which is denoted by $\textrm{Max}(G)$. In this article, we obtain upper bounds for $\textrm{Max}(G)$, for $G=G_{n,m}$, grids with $n,m \geq 9$, $m\equiv 1 \mod 4$ and $nm$ odd.
\end{abstract}
%
\maketitle
\section{Introduction}
Let $G$ be a graph with edge and vertex sets $E(G)$ and $V(G)$, respectively, and for any $u,v \in V(G)$, let $d(u,v)$ be the distance between $u$ and $v$. A matching is a set of  edges with no shared vertices.  The vertices incident to an edge of a
matching $M$ are said to be saturated by $M$, the other vertices are unsaturated by M.
  A subset $M \subseteq E(G)$ of $G$ is an induced matching of $G$ if for any two edges $e_1=u_iu_j$ and $e_2=v_iv_j$ in $M$, then $d(u_i,v_i) \geq d(u_i,v_j) \geq 2$ and $d(u_j,v_i) \geq d(u_j,v_j) \geq 2$. In other words, $M$ is an induced matching of $G$ if for any two edges $e_1,e_2$ in $M$, there is no edge in $G$ connecting $e_1$ and $e_2$. Equivalently, an induced matching is a matching which forms an induced subgraph.

Introduced by Stockmeyer and Vazirani \cite{SV1} as a special case of the well known matching problem, the concept finds applications, among others, in cryptology where certain communication channels between two ends are classified \cite{CST1}.
	
	The size $|M|$ of an induced matching $M$ is the number of edges in the induced matching. Denote by $\textrm{Max}(G)$ the maximum size of an induced matching in $G$.  A maximum induced matching $M$ in $G$ is an induced matching with $\textrm{Max}(G)$ edges and $\textrm{Max}(G)$ is known as  maximum induced matching number (or strong matching number) of $G$. 
 Unlike in the case with finding the maximum matching number of a graph, which can be obtained in polynomial time \cite{E1}, obtaining $\textrm{Max}(G)$ in general, is $NP-$hard, even for classes of graph such as the regular bipartite graph \cite{DD1}.
	
	In \cite{Z1}, it was observed that for any $G$ with maximum degree $\Delta (G)$,
	
	\begin{equation}
	\textrm{Max}(G) \geq \frac{|V(G)|}{2(2\Delta (G)^2 +2\Delta (G)+1)}
	\end{equation}
	The above bound is certainly not a sharp one and therefore, Joos in \cite{J1}, presented the following bound:
	
	\begin{equation}
	\textrm{Max}(G) \geq \frac{|V(G)|}{(\lceil \frac{\Delta(G)}{2}\rceil + 1)(\lfloor \frac{\Delta(G)}{2}\rfloor +1)}
	\end{equation}
	The last bound is subject to the size of $\Delta (G)$. He also showed that it holds for $\Delta(G) = 1000$ and this, according to him, could be reviewed down to $200$. In the end, he conjectured that this bound holds for $\Delta(G) \geq 3$ with the exception of certain graphs that he listed. Inspired by the work in \cite{J1}, Nguyen \cite{N1} showed that the conjecture is true for $\Delta(G)=4$ as long as $G$ is not one of the exception graphs.
	
	The maximum induced matching number of many graphs can be obtained efficiently just as in the cases of chordal graphs \cite{CST1}, bounds of bounded cliques width, intersection graphs \cite{C1}, circular arc graphs \cite{GL1} among others.
	
 	A grid $G_{n,m}$ is obtained by the Cartesian product of any two paths of lengths $n$ and $m$, where $n, m \geq 2$ are integers, representing rows and columns of the grid, respectively. We introduce an odd grid as a grid whose path factors are of odd order. Marinescu-Ghemaci in \cite{RMG1}, obtained $\textrm{Max}(G)$ values for all grids with even $nm$, and some cases where $nm$ is odd. She also gave useful lower and upper bounds. Particularly, she showed that for any odd grid $G_{n,m}$, $\textrm{Max}(G_{n,m}) \leq \left\lfloor \frac{nm+1}{4} \right\rfloor$.

	This paper improves Marinescu-Ghemachi's upper bound for $G_{n, m}$, $n, m$ odd, $m \equiv 1 \mod 4$. 
The results provide new upper bounds for some cases whose lower bounds were established in \cite{RMG1} and thus, in a number of situations, precise values of $\textrm{Max}(G_{n, m})$ were obtained. These may also prove useful in probing some of the conjectures left in \cite{RMG1}.
	
	\section{Definitions and Preliminary Results}
	Grid, $G_{n,m}$, as defined in this work, is the Cartesian product of paths $P_n$ and $P_m$ with $n$ and $m$ being positive integers, where $P_n$ and $P_m$ have disjoint vertex sets $V(P_n)=\left\{u_1,u_2,\cdots, u_n\right\}$ and $V(P_m)=\left\{v_1,v_2,\cdots, v_m\right\} $, respectively. Unless explicitly stated, $n$ is any odd integer, $m\equiv 1\mod 4$ and $2\leq n \leq m$. We introduce the following notations: $V_i=\left\{u_1v_i, u_2v_i, \cdots, u_nv_i\right\} \subset V(G_{n,m})$ and $U_i=\left\{u_iv_1, u_iv_2, \cdots, u_iv_m\right\}\subset V(G_{n,m})$; for edge set $E(G_{n,m})$ of $G_{n,m}$, if $u_iv_j \; u_kv_j \in E(G_{n,m})$ and $u_iv_j \; u_iv_k \in E(G_{n,m})$, we write $u_{\left\{i,k\right\}}v_j \in E(G_{n,m})$ and $u_iv_{\left\{j,k\right\}} \in E(G_{n,m})$ respectively.
	
	A vertex $v$ is said to be saturated by an induced matching $M$ if it is a member of an edge in $M$ and unsaturated by $M$, otherwise. We say $v$ is saturable if $v$ is saturated by $M$ or if $v$ is unsaturated by $M$, then for every saturated vertex $u$ in $V(G)$, $d(v,u) \geq 2$. This implies that an unsaturated vertex can become saturated if there is no saturated vertex within distance $2$.  If $v$ can not be saturated, then we say $v$ is an isolated saturable vertex. A boundary vertex is a vertex on any of $V_1,V_n,U_1, U_m$. As we shall see in this work however, boundaries of some subgraph are defined. A saturable vertex $v$ in subgraph $G_a$ of $G$, which is not saturated by induced matching $M_a$ of $G_a$ can still be saturated by $M$ of $G$ in case $v$ is on the boundary of $G_a$.   However, if $v$ is not on the boundary of $G_a$, then $v$ is isolated. The sets of all saturated vertices, saturable vertices and isolated vertices in a graph $G$ are denoted by $V_s(G),V_{sb}(G)$ and $V_{is}(G)$, respectively. Clearly, $\left|V_s(G)\right|$ is even and $V_s(G) \subseteq V_{sb}(G)$ .
	{\tiny{
\begin{center}
\pgfdeclarelayer{nodelayer}
\pgfdeclarelayer{edgelayer}
\pgfsetlayers{nodelayer,edgelayer}
\begin{tikzpicture}
	\begin{pgfonlayer}{nodelayer}

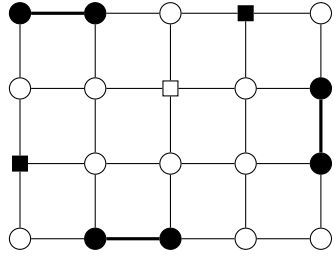
	\node [minimum size=0cm,]  at (3,0.5) {Saturable vertices as black squares and an isolated vertex as white square in an induced matching of $G_{4,5}$};

		\node [minimum size=0cm,draw,circle] (1) at (1,1) {};
		\node [minimum size=0cm,draw,fill=black!, circle] (2) at (2,1) {};
		\node [minimum size=0cm, draw,fill=black!, circle] (3) at (3,1) {};
		\node [minimum size=0cm,draw,circle] (4) at (4,1) {};
		\node [minimum size=0cm,draw,circle] (5) at (5,1) {};

		\node [minimum size=0cm, draw,fill=black!, rectangle] (6) at (1,2) {};
		\node [minimum size=0cm,draw,circle] (7) at (2,2) {};
		\node [minimum size=0cm,draw,circle] (8) at (3,2) {};
		\node [minimum size=0cm,draw,circle] (9) at (4,2) {};
		\node [minimum size=0cm,draw,fill=black!, circle] (10) at (5,2) {};

		\node [minimum size=0cm,draw,circle] (11) at (1,3) {};
		\node [minimum size=0cm,draw,circle] (12) at (2,3) {};
		\node [minimum size=0cm,draw,rectangle] (13) at (3,3) {};
		\node [minimum size=0cm,draw,circle] (14) at (4,3) {};
		\node [minimum size=0cm,draw,fill=black!, circle] (15) at (5,3) {};

		\node [minimum size=0cm,draw,fill=black!, circle] (16) at (1,4) {};
		\node [minimum size=0cm,draw,fill=black!, circle] (17) at (2,4) {};
	  \node [minimum size=0cm,draw,circle] (18) at (3,4) {};
		\node [minimum size=0cm,draw,fill=black!,rectangle] (19) at (4,4) {};
		\node [minimum size=0cm,draw,circle] (20) at (5,4) {};

		\end{pgfonlayer}
	\begin{pgfonlayer}{edgelayer}
	\draw [thin=1.00] (1) to (2);
		\draw [very thick=1.00] (2) to (3);
		\draw [thin=1.00] (3) to (4);
		\draw [thin=1.00] (4) to (5);
		
		\draw [thin=1.00] (6) to (7);
		\draw [thin=1.00] (7) to (8);
		\draw [thin=1.00] (8) to (9);
		\draw [thin=1.00] (9) to (10);
		
		\draw [thin=1.00] (11) to (12);
		\draw [thin=1.00] (12) to (13);
		\draw [thin=1.00] (13) to (14);
		\draw [thin=1.00] (14) to (15);
		
		\draw [very thick=1.00] (16) to (17);
		\draw [thin=1.00] (17) to (18);
		\draw [thin=1.00] (18) to (19);
		\draw [thin=1.00] (19) to (20);

		\draw [thin=1.00] (1) to (6);
		\draw [thin=1.00] (6) to (11);
		\draw [thin=1.00] (11) to (16);
		
		\draw [thin=1.00] (2) to (7);
		\draw [thin=1.00] (7) to (12);
		\draw [thin=1.00] (12) to (17);
		
		\draw [thin=1.00] (3) to (8);
		\draw [thin=1.00] (8) to (13);
		\draw [thin=1.00] (13) to (18);
		
		\draw [thin=1.00] (4) to (9);
		\draw [thin=1.00] (9) to (14);
		\draw [thin=1.00] (14) to (19);
		
		\draw [thin=1.00] (5) to (10);
		\draw [very thick=1.00] (10) to (15);
		\draw [thin=1.00] (15) to (20);

	\end{pgfonlayer}
\end{tikzpicture}
\end{center}
}}

The following results about grid $G_{n,m}$ are from \cite{RMG1}:
\begin{lem}\label{lem2.1} Let $m,n \geq 2$ be two positive integers.
\begin{enumerate}
\item
If $m \equiv 2 \mod 4$ and $n$ odd, then, $|V_{sb}(G_{n,m})|=\frac{mn+2}{2}$ and $|V_{sb}(G_{n,m})|=\frac{mn}{2}$ otherwise;
\item For $m\geq 3$, $m$ odd, $|V_{sb}(G_{n,m})|=\frac{nm+1}{2}$, for $n\in \left\{3,5\right\}$.
\end{enumerate}
\end{lem}
\begin{thm}\label{thm2.2} For $G_{n,m}$ where $2\leq n \leq m$,  let $|M|= \textrm{Max}(G_{n,m})$. Then, for $n$ even, $\textrm{Max}(G_{n,m})=\left\lceil \frac{mn}{4}\right\rceil$, for $n \in \left\{3,5\right\}$, $m \equiv 1 \mod 4$, $\textrm{Max}(G_{n,m})=\frac{n(m-1)}{4}+1$ and for $m \equiv 3 \mod 4$, $\textrm{Max}(G_{n,m})=\frac{n(m-1)+2}{4}$.
\end{thm}	
\begin{rem}\label{rem2.3} For $m \equiv 1 \mod 4$, $|V_{sb}(G_{3,m})|=$ $2(\textrm{Max }(G_{3,m}))=|V_s(G_{n,m})|$.
\end{rem}
\begin{thm}\label{thm2.4} For $m, n$ odd integers, $\textrm{Max}(G_{n,m}) \leq \left\lfloor \frac{mn+1}{4}\right\rfloor$.
\end{thm}
The obvious implication of Theorem \ref{thm2.4}, based on the proof, is that $|V_{sb}G_{n,m}| \leq \frac{nm+1}{2}$, for $n, m$ odd.

\section{Results}
We start our results by stating a few observations.
 \begin{rem}\label{rem3.1} Let $G_{3,3}$ be a $3 \times 3$ grid with induced matching $M$. Clearly, by Lemma \ref{lem2.1}, $|V_{sb}(G_{3,3})|=5$. Suppose $u_{\{1,2\}}v_2 \in M$, then $|M|=1$. However, there are non adjacent saturable vertices $u_3v_1$ and $u_3v_3$.
 \end{rem}

\begin{lem}\label{lem3.2} Let $M$ be an induced matching of $G_{3,m}$. If $u_{\{1,2\}}v_2 \in M$. Then, $|M| \neq \textrm{Max}(G_{3,m})$.
\end{lem}

\begin{proof}  Suppose $u_{\{1,2\}}v_2 \in M$ and let $G_a=G_{3,m-3} \subset G_{3,m}$, be a subgrid of $G_{3,m}$ induced by $V(G_{3,m}) \backslash \left\{V_1,V_2,V_3\right\},$ where $V_1,V_2,V_3 \subset V(G_{3,m})$, $m \geq 3$.

 \textbf{Case I}: Suppose $m \equiv 3 \mod 4$, $m \geq 7$. Then $m=4k+3$ for some positive integers $k$. By Lemma \ref{lem2.1} and Theorem \ref{thm2.2}, $|V_s(G_a)|=6k$. Now, with $u_{\{1,2\}}v_2 \in M$, $u_3v_1$ is saturable, by its position and $V_s(G_a)=V_{sb}(G_a)$.  Since $m-3$ is even, then either $u_3v_3$ remains isolated (or unsaturated) or if it is forced to be saturated with $u_3v_4$, a saturable vertex in $V(G_a)$ becomes isolated (or unsaturated). So without loss of generality, we set $u_3v_3$ as unsaturated. Therefore $G_{3,m}$ contains $3+6k$ saturable vertices and then, $|M| \leq 3k+1$, which is a contradiction since by Theorem  \ref{thm2.2}, $\textrm{Max}(G_{3,m})=3k+2$. Note that $m=3$, has been covered by Remark \ref{rem3.1}.

 \textbf{Case II}: Suppose $m \equiv 1 \mod 4$. From Lemma \ref{lem2.1}, Theorem \ref{thm2.2} and following similar argument as in Case I, we see that if $u_{\{1,2\}}v_2 \in M$ then one vertex in $V_{sb}(G_3,m)$ will become isolated and therefore, $|V_{sb}(G_{3,m})|=6k+1$. Hence $|M| \leq 3k$ and a contradiction since $\textrm{Max} (G_{3,m})=3k+1$ if $m=4k+1$.

 \textbf{Case III}: If $m \equiv 0 \mod 4$, then $m-3 \equiv 1 \mod 4$, and by Remark \ref{rem2.3}, $G_{3,m-3}=G_a$, $|V_s(G_a)|=|V_{sb}(G_a)|$. Now $|V_s(G_a)|= \frac{3(m-3)+1}{2}=6k-4$. By following the argument in the previous cases, we see that $|V_{sb}(G_a)|=6k-1$. Therefore, $|M| \leq 3k-1$, which is less than $3k$.

\textbf{Case IV}: If $m \equiv 2 \mod 4$, then $m-3 \equiv 3 \mod 4$ and therefore, $|V_{sb}(G_a)|=\frac{3(4k+2-3)+1}{2}=6k-1$, which is odd. Therefore, there exists a saturable vertex in $V(G_a)$, which can pair with $u_3v_3$ and thus the two vertices become saturated. This way, we have the total saturable vertices in $G_{3,m}$, which has $u_{\{1,2\}}v_2 \in M$, to be $6k+3$ which implies that $|M| \leq 3k+1 $. But from the results in Theorem \ref{thm2.2}, if $m \equiv 2 \mod 4$, $\textrm{Max}(G_{3,m})=3k+2$. This is a contradiction and the claim holds.
\end{proof}

\begin{rem}\label{rem3.2} Suppose we have $u_{\{1,2\}}v_2 \in M$ in a grid $G_{3,m}$,  then the following holds for $|V_{sb}(G_{3,m})|$, from Lemma \ref{lem3.2}.
\begin{center}

 $ \begin{array}{|c|c|} \hline
 m & |V_{sb}(G_{3,m})|\\ \hline
 4k & 6k-1 \\ 
 4k+1 & 6k+1 \\ 
 4k+2 &6k+3  \\ 
 4k+3 & 6k+3 \\ \hline

 \end{array}$

 \end{center}

\end{rem}
\begin{lem}\label{lem3.3} If $m\equiv 1 \mod 4$ and there exists $M$, an induced matching of $G_{3,m}$, such that $u_{\{1,2\}}v_j, u_{\{1,2\}}v_{j+2} \in M$, then $|M| \neq \textrm{Max}(G_{3,m})$.
\end{lem}

\begin{proof} Let $u_{\{1,2\}}v_j, u_{\{1,2\}}v_{j+2} \in M$, where $M$ is an induced matching of $G_{3,m}$.

\textbf{Case I}: Suppose that $j+1 \equiv 3 \mod 4$, then $m-(j+1) \equiv 2 \mod 4$. Since $u_{\{1,2\}}v_j \in M$, by Lemma \ref{lem3.2}, suppose there exist an induced matching $M'$ in $G_a=G_{3,j+1}$, induced by $V_1, V_2,\cdots, V_{j+1}$, with $u_{\{1,2\}}v_j \in M'$, then $|M'| \neq \textrm{Max}(G_a)$. By Remark \ref{rem3.2}, given a non-negative integer $l$, $|V_{sb}(G_a)|= 6l+3$, which being odd, contains a saturable vertex $v'=u_3v_{j+1}$ which is not a member of $V_s(G_a)$. In fact, $v'\in V_{is}(G_a)$, since $u_{\{1,2\}}v_{j+1} \in M$.  Since $m-(j+1) \equiv 2 \mod 4$, then given a subgrid $G_b=G_{3,m-(j+1)}$, induced by $V_{j+2}, V_{j+3}, \cdots , V_m$ $|V_s(G_b)|=6(k-l)-2$. Therefore, $|V_{sb}(G_{3,m})|=6k$ and hence, $|M| \leq 3k$, which is less that $3k+1$.

\textbf{Case II}: If $j+1 \equiv 1 \mod 4$, by following the argument in Case I, $v' \in V_{is} (G_{3,m})$. Now, $j+1=4l+1$, and by the isolation of $v'$, and Remark \ref{rem3.2}, $|V_{sb}(G_{3,m})|=6l+|V_{sb}(G_b)|$. Meanwhile, $m-(j+1) \equiv 0 \mod 4$ and therefore, $|V_{sb}(G_b)|=6k$. Thus, $|M| \leq 3k$.

\textbf{Case III}: If $j+1$ is even, we follow similar arguments as the earlier cases.

\end{proof}

\begin{lem} Suppose that $m \equiv 1 \mod 4$ and $M$ is an induced matching of $G_{3,m}$. If $u_{\{1,2\}}u_j, u_{\{1,2\}}u_{j+3} \in M$, then $|M|\neq \textrm{Max}(G_{3,m})$.
\end{lem}

\begin{proof} For some positive integers $k$, let $m=4k+1$. Now suppose that $j \equiv 3 \mod 4$. This implies that $j+3 \equiv 2 \mod 4$ and $u_{\{1,2\}}u_j, u_{\{1,2\}}u_{j+3} \in M$. Let $G_a=G_{3,j+1}$, $G_b=G_{3,m-(j+1)} \subset G_{3,m}$, induced by $V_1,V_2,\cdots, V_{j+1}$ and $V_{j+2}, V_{j+3}, \cdots, V_m$, respectively. By earlier result and remark, $|V_{sb}(G_a)|=6l-1$ since $j+1=4l$. Also, $|V_{sb}(G_b)|=6(k-l)+1$, since $m-(j+1)=4(k-l)+1$. Certainly, $|V_{sb}(G_{3,m})|=6k$. Thus $|M| \leq 3k$. Therefore, $|M| \neq 3k$. For $j \equiv 1 \mod 4$, $j+3 \equiv 3 \mod 4$. Let $u_{\{1,2\}}u_j, u_{\{1,2\}}u_{j+3} \in M$. By earlier lemma and result, we have that $|V_{sb}(G_a)|=6l+3$. Since $j+1=4l+2$. Also, $|V_{sb}(G_b)|=6[(k-l)-1]+3$. Therefore, $|V_{sb}(G_{3,m})|=6k$ and therefore, $|M|=3k$, which is a contradiction.

\end{proof}

\begin{rem}\label{rem3.3} By following similar argument as in the last result, it is easy to see that if $M$ contains $u_{\{1,2\}}u_j, u_{\{1,2\}}u_{j+4}$, then $|M|\neq \textrm{Max} (G_{3,m})$. Therefore, suppose $M'$ is a maximum induced matching of $(G_{3,m})$ and  $u_{\{1,2\}}u_j, u_{\{1,2\}}u_{j+8} \in M'$, then there is no $u_{\{1,2\}}u_{j+k} \in M$, such that $2 \leq k \leq 6$.
\end{rem}

\begin{thm}\label{thm3.4} For $m=4k+1$, there exist at least $2k$ saturated vertices in $U_1 \subset V(G_{3,m})$.
\end{thm}

\begin{proof} Let $G_a=G_{2,m} \subset G_{3,m}$, induced by $U_2,U_3$, a subgrid of $G_{3,m}$. From earlier results, $|V_{sb}(G_a)|=4k+2$. Now, $|V_{sb}(G_{3,m})|=6k+2$. Therefore, $V_1$ has at least $2k$ saturable vertices.
\end{proof}

\begin{thm}\label{thm3.5} Suppose $u_{\{1,2\}}u_j,u_{\{1,2\}}v_{j+8} \in M$, in  a $G_{3,m}$ and $m \equiv 1 \mod 4$.
\begin{enumerate}
\item[(a)] There exists four other saturated vertices from $u_1v_j$ to $u_1v_{j+8}$.
\item[(b)] There exists at most one saturated vertex between $u_2v_j$ and $u_2v_{j+8}$.
\item[(c)] There exists at most five saturated vertices from $v_3v_j$ to $u_3v_{j+8}$.
\end{enumerate}

\end{thm}
 \begin{proof} \begin{enumerate} \item[(a)] From Remark \ref{rem3.3}, let $|M|$ = $\textrm{Max}(G_{3,m})$, and $u_{\{1,2\}}u_j,u_{\{1,2\}}v_{j+8} \in M$, then there is no $u_{\{1,2\}}v_k \in E(G_{3,m})$, $1 \leq k \leq 7$ such that $u_{\{1,2\}}v_k \in M$. Thus, suppose there exists another saturated vertex such that 
 $u_1v_{\{j+k,j+k-1\}} \in M$. Then there exist at least four saturated vertices from $u_1v_1$ to $u_1v_m$. Now suppose that there is no other saturated vertex on $U_1$, then by the Theorem \ref{thm3.4}, the grid $G_a=G_{2,m} \subset G_{3,m}$, induced by vertices $\left\{u_2v_j, u_2v_{j+1},  \cdots, u_2v_{j+8}\right\}$ and $\left\{u_3v_j, u_3v_{j+1},  \cdots, u_3v_{j+8}\right\}$ must contain ten saturated vertices, (including $u_2v_j$ and $u_2v_{j+8}$). Clearly, vertices $u_2v_{j+1},u_2v_{j+7}, u_3v_j$ and $u_3v_{j+8}$ can not be saturated by $G_a$. It is clear, therefore, that $G_a$ only has 8 saturable vertices, which is a contradiction. Thus, there exists two more saturated vertices in $U_1$, and hence the claim.
\item[ ] Parts $(b)$ and $(c)$ follow from $(a)$.
 \end{enumerate}
 \end{proof}

 \begin{rem}\label{rem3.4} \begin{enumerate} \item[(a)] Since there are five saturated vertices between $u_3v_j$ and $u_3v_{j+8}$ and there exist only one saturated vertex between $u_2v_j$ and $u_2v_{j+8}$, then there are edges $e_1, e_2 \in E(G_b)$, $G_b$ induced by $u_3v_j, u_3v_{j+1}, \cdots u_3v_{j+8}$.
 \item[(b)] Suppose $m= 4k+1$, $k$ being positive integers, then, at $G_c=G_{1,m} \subset G_{3,m}$, induced by either $V_1$ or $V_3$, there exists at least $k$ edges of $G_c$ in the maximum induced matching $M$ of $(G_{3,m})$.
     \end{enumerate}
 \end{rem}

 {\tiny{
\begin{center}
\pgfdeclarelayer{nodelayer}
\pgfdeclarelayer{edgelayer}
\pgfsetlayers{nodelayer,edgelayer}
\begin{tikzpicture}
	\begin{pgfonlayer}{nodelayer}
	
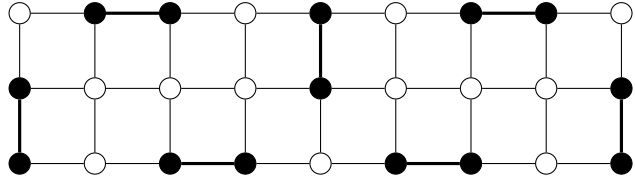
	\node [minimum size=0cm,]  at (-9,6.5) { A $G_{3,9}$ grid with $Max(G_{3,9})=7$};

		\node [minimum size=0cm,fill=black!,draw,circle] (1) at (-13,7) {};

		\node [minimum size=0cm,draw,circle] (2) at (-12,7) {};
		\node [minimum size=0cm,fill=black!,draw,circle] (3) at (-11,7) {};
		\node [minimum size=0cm,fill=black!,draw,circle] (4) at (-10,7) {};
		\node [minimum size=0cm,draw,circle] (5) at (-9,7) {};
		\node [minimum size=0cm,draw,fill=black!,circle] (6) at (-8,7) {};
		\node [minimum size=0cm,draw,fill=black!,circle] (7) at (-7,7) {};
		\node [minimum size=0cm,draw,circle] (8) at (-6,7) {};
		\node [minimum size=0cm,draw,fill=black!,circle] (9) at (-5,7) {};
		
		\node [minimum size=0cm,draw,fill=black!,circle] (11) at (-13,8) {};
		\node [minimum size=0cm,draw,circle] (12) at (-12,8) {};
		\node [minimum size=0cm,draw,circle] (13) at (-11,8) {};
		\node [minimum size=0cm,draw,circle] (14) at (-10,8) {};
		\node [minimum size=0cm,draw,fill=black!,circle] (15) at (-9,8) {};
		\node [minimum size=0cm,draw,circle] (16) at (-8,8) {};
		\node [minimum size=0cm,draw,circle] (17) at (-7,8) {};
		\node [minimum size=0cm,draw,circle] (18) at (-6,8) {};
		\node [minimum size=0cm,draw,fill=black!,circle] (19) at (-5,8) {};
		
	  \node [minimum size=0cm,draw,circle] (21) at (-13,9) {};
		\node [minimum size=0cm,draw,fill=black!,circle] (22) at (-12,9) {};
		\node [minimum size=0cm,draw,fill=black!,circle] (23) at (-11,9) {};
		\node [minimum size=0cm,draw,circle] (24) at (-10,9) {};
		\node [minimum size=0cm,draw,fill=black!,circle] (25) at (-9,9) {};
		\node [minimum size=0cm,draw,circle] (26) at (-8,9) {};
		\node [minimum size=0cm,draw,fill=black!,circle] (27) at (-7,9) {};
		\node [minimum size=0cm,draw,fill=black!,circle] (28) at (-6,9) {};
		\node [minimum size=0cm,draw,circle] (29) at (-5,9) {};

		\end{pgfonlayer}
	\begin{pgfonlayer}{edgelayer}
		\draw [thin=1.00] (1) to (2);
		\draw [thin=1.00] (2) to (3);
		\draw [very thick=1.00] (3) to (4);
		\draw [thin=1.00] (4) to (5);
		\draw [thin=1.00] (5) to (6);
		\draw [very thick=1.00] (6) to (7);
		\draw [thin=1.00] (7) to (8);
		\draw [thin=1.00] (8) to (9);
		\draw [thin=1.00] (11) to (12);
		\draw [thin=1.00] (12) to (13);
		\draw [thin=1.00] (13) to (14);
		\draw [thin=1.00] (14) to (15);
		\draw [thin=1.00] (15) to (16);
		\draw [thin=1.00] (16) to (17);
		\draw [thin=1.00] (17) to (18);
		\draw [thin=1.00] (18) to (19);
		\draw [thin=1.00] (21) to (22);
		\draw [very thick=1.00] (22) to (23);
		\draw [thin=1.00] (23) to (24);
		\draw [thin=1.00] (24) to (25);
		\draw [thin=1.00] (25) to (26);
		\draw [thin=1.00] (26) to (27);
		\draw [very thick=1.00] (27) to (28);
		\draw [thin=1.00] (28) to (29);

		\draw [very thick=1.00] (1) to (11);
		\draw [thin=1.00] (11) to (21);
		
		\draw [thin=1.00] (2) to (12);
		\draw [thin=1.00] (12) to (22);

		\draw [thin=1.00] (3) to (13);
		\draw [thin=1.00] (13) to (23);

		\draw [thin=1.00] (4) to (14);
		\draw [thin=1.00] (14) to (24);

		\draw [thin=1.00] (5) to (15);
		\draw [very thick=1.00] (15) to (25);

		\draw [thin=1.00] (6) to  (16);
		\draw [thin=1.00] (16) to (26);

		\draw [thin=1.00] (7) to  (17);
		\draw [thin=1.00] (17) to (27);

		\draw [thin=1.00] (8) to  (18);
		\draw [thin=1.00] (18) to (28);

		\draw [very thick=1.00] (9) to (19);
		\draw [thin=1.00] (19) to (29);

	\end{pgfonlayer}
\end{tikzpicture}
\end{center}
}}

 Next, we consider the grid $G_{4,m}$, $m \equiv 1 \mod 4$.

 \begin{lem}\label{lem3.6} Let $|M|= \textrm{Max} (G_{4,m})$ and let $U_4$ contain $\frac{m-1}{2}$ saturated vertices. Then, for any edge $e_1 \in E(G_a)$, $G_a=(G_{1,m}) \subset G_{4,m}$, induced by $U_4 \subset V(G_{4,m})$, $e_1 \notin M$.
 \end{lem}

 \begin{proof} Let $u_4v_i, u_4v_{i+1} \in U_4 \subset V(G_{4,m})$, be saturated vertices. By hypothesis, there are $\frac{m-5}{2}$ other saturated vertices in $U_4$. Suppose $G_a=G_{2,m} \subset G_{4,m}$, induced by $U_1, U_2$. The $|V_s(G_a)|=|V_{sb}(G_a)|=m+1$. Let $G_b=G_{2,m} \subset G_{4,m}$, induced by $U_3,U_4$. Since $|V_s(G_a)|=m+1$, then $|V_s(G_b)|\leq m-1$  where $|V_s(G_{4,m})|=2m$. By hypothesis, suppose there exist $\frac{m-1}{2}$ saturated vertices in $U_4$, then there are also at most $\frac{m-1}{2}$ saturated vertices in $U_3$. Without loss of generality, suppose for all the $\frac{m-5}{2}$ other saturated vertices in $U_s$, there exist adjacent saturated vertices in $U_3$, then we have that $|V_s(G_b)| \geq m-3$.

 \textbf{Claim:} There are at most $\frac{m-5}{2}$ saturated vertices in $U_3$.

 \textbf{Reason:} Since there is at most $\frac{m-1}{2}$ saturated vertices in $U_4$, then suppose $v_k$ is saturable in $U_3$, $v_k$ is not incident to a saturable vertex in $U_4$. Also, since $V_s(G_a)=V_{sb}(G_a)$, then there is no saturable vertex in $U_2$ to which $v_k$ is incident to form an edge in $M$. Thus, suppose there exists a vertex $v_{k-1}$, adjacent to $v_k \in U_3$, $v_{k-1}$ will be adjacent to a saturated vertex in $U_2$ since there can not be two adjacent vertices both of which are not saturated in $U_2$. Thus, $v_k$ is isolated. \\
 Finally, $|V_s(G_{4,m})|\leq 2m-2$ which implies that $|M| \leq m-1$. However by Theorem \ref{thm2.2}, $\textrm{Max}(G_{4,m})=m$.
 \end{proof}

\begin{corollary}\label{cor3.7} Let $G_{n,m}$ be a grid with $n \equiv 0 \mod 4, m \equiv 1 \mod 4$ and $U_n$ contains $\frac{m-1}{2}$ saturable vertices, with $|M|= \textrm{Max}({G}_{n,m})$, then no two saturated vertices, say, $v',v'' \in U_n$ such that $v'v'' \in M$.
\end{corollary}

{\tiny{
\begin{center}
\pgfdeclarelayer{nodelayer}
\pgfdeclarelayer{edgelayer}
\pgfsetlayers{nodelayer,edgelayer}
\begin{tikzpicture}
	\begin{pgfonlayer}{nodelayer}
	
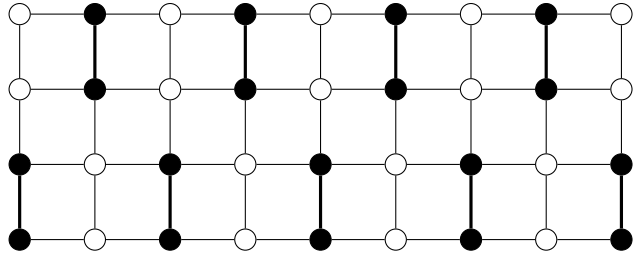
	\node [minimum size=0cm,]  at (-9,6.5) { A $G_{4,9}$ Grid with $Max(G_{4,9})=11$};

		\node [minimum size=0cm,fill=black!,draw,circle] (1) at (-13,7) {};
		\node [minimum size=0cm,draw,circle] (2) at (-12,7) {};
		\node [minimum size=0cm,fill=black!,draw,circle] (3) at (-11,7) {};
		\node [minimum size=0cm,draw,circle] (4) at (-10,7) {};
		\node [minimum size=0cm,fill=black!,draw,circle] (5) at (-9,7) {};
		\node [minimum size=0cm,draw,circle] (6) at (-8,7) {};
		\node [minimum size=0cm,fill=black!,draw,circle] (7) at (-7,7) {};
		\node [minimum size=0cm,draw,circle] (8) at (-6,7) {};
		\node [minimum size=0cm,fill=black!,draw,circle] (9) at (-5,7) {};
		\node [minimum size=0cm,fill=black!,draw,circle] (11) at (-13,8) {};
		\node [minimum size=0cm,draw,circle] (12) at (-12,8) {};
		\node [minimum size=0cm,fill=black!,draw,circle] (13) at (-11,8) {};
		\node [minimum size=0cm,draw,circle] (14) at (-10,8) {};
		\node [minimum size=0cm,fill=black!,draw,circle] (15) at (-9,8) {};
		\node [minimum size=0cm,draw,circle] (16) at (-8,8) {};
		\node [minimum size=0cm,fill=black!,draw,circle] (17) at (-7,8) {};
		\node [minimum size=0cm,draw,circle] (18) at (-6,8) {};
		\node [minimum size=0cm,fill=black!,draw,circle] (19) at (-5,8) {};
	  \node [minimum size=0cm,draw,circle] (21) at (-13,9) {};
		\node [minimum size=0cm,fill=black!,draw,circle] (22) at (-12,9) {};
		\node [minimum size=0cm,draw,circle] (23) at (-11,9) {};
		\node [minimum size=0cm,fill=black!,draw,circle] (24) at (-10,9) {};
		\node [minimum size=0cm,draw,circle] (25) at (-9,9) {};
		\node [minimum size=0cm,fill=black!,draw,circle] (26) at (-8,9) {};
		\node [minimum size=0cm,draw,circle] (27) at (-7,9) {};
		\node [minimum size=0cm,fill=black!,draw,circle] (28) at (-6,9) {};
		\node [minimum size=0cm,draw,circle] (29) at (-5,9) {};
		\node [minimum size=0cm,draw,circle] (31) at (-13,10) {};
		\node [minimum size=0cm,fill=black!,draw,circle] (32) at (-12,10) {};
		\node [minimum size=0cm,draw,circle] (33) at (-11,10) {};
		\node [minimum size=0cm,fill=black!,draw,circle] (34) at (-10,10) {};
		\node [minimum size=0cm,draw,circle] (35) at (-9,10) {};
		\node [minimum size=0cm,fill=black!,draw,circle] (36) at (-8,10) {};
		\node [minimum size=0cm,draw,circle] (37) at (-7,10) {};
		\node [minimum size=0cm,fill=black!,draw,circle] (38) at (-6,10) {};
		\node [minimum size=0cm,draw,circle] (39) at (-5,10) {};

		\end{pgfonlayer}
	\begin{pgfonlayer}{edgelayer}
		\draw [thin=1.00] (1) to (2);
		\draw [thin=1.00] (2) to (3);
		\draw [thin=1.00] (3) to (4);
		\draw [thin=1.00] (4) to (5);
		\draw [thin=1.00] (5) to (6);
		\draw [thin=1.00] (6) to (7);
		\draw [thin=1.00] (7) to (8);
		\draw [thin=1.00] (8) to (9);
		\draw [thin=1.00] (11) to (12);
		\draw [thin=1.00] (12) to (13);
		\draw [thin=1.00] (13) to (14);
		\draw [thin=1.00] (14) to (15);
		\draw [thin=1.00] (15) to (16);
		\draw [thin=1.00] (16) to (17);
		\draw [thin=1.00] (17) to (18);
		\draw [thin=1.00] (18) to (19);
		\draw [thin=1.00] (21) to (22);
		\draw [thin=1.00] (22) to (23);
		\draw [thin=1.00] (23) to (24);
		\draw [thin=1.00] (24) to (25);
		\draw [thin=1.00] (25) to (26);
		\draw [thin=1.00] (26) to (27);
		\draw [thin=1.00] (27) to (28);
		\draw [thin=1.00] (28) to (29);
		
		\draw [thin=1.00] (31) to (32);
		\draw [thin=1.00] (32) to (33);
		\draw [thin=1.00] (33) to (34);
		\draw [thin=1.00] (34) to (35);
		\draw [thin=1.00] (35) to (36);
		\draw [thin=1.00] (36) to (37);
		\draw [thin=1.00] (37) to (38);
		\draw [thin=1.00] (38) to (39);

		\draw [very thick=1.00] (1) to (11);
		\draw [thin=1.00] (11) to (21);
		\draw [thin=1.00] (21) to (31);

		\draw [thin=1.00] (2) to (12);
		\draw [thin=1.00] (12) to(22);
		\draw [very thick=1.00] (22) to (32);

		\draw [very thick=1.00] (3) to (13);
		\draw [thin=1.00] (13) to (23);
		\draw [thin=1.00] (23) to (33);
		
		\draw [thin=1.00] (4) to (14);
		\draw [thin=1.00] (14) to (24);
		\draw [very thick=1.00] (24) to (34);
		
		\draw [very thick=1.00] (5) to (15);
		\draw [thin=1.00] (15) to (25);
		\draw [thin=1.00] (25) to (35);
		
		\draw [thin=1.00] (6) to  (16);
		\draw [thin=1.00] (16) to (26);
		\draw [very thick=1.00] (26) to (36);

		\draw [very thick=1.00] (7) to  (17);
		\draw [thin=1.00] (17) to (27);
		\draw [thin=1.00] (27) to (37);
		
		\draw [thin=1.00] (8) to  (18);
		\draw [thin=1.00] (18) to (28);
		\draw [very thick=1.00] (28) to (38);
		
		\draw [very thick=1.00] (9) to (19);
		\draw [thin=1.00] (19) to (29);
		\draw [thin=1.00] (29) to (39);

	\end{pgfonlayer}
\end{tikzpicture}
\end{center}
}}

Next we observe some results in the grid $G_{5,m}$, where $m \equiv 1 \mod 4 $.

\begin{lem}\label{lem3.8} Let $G_{5,5}$ be a grid and suppose that $|M|=\textrm{Max} (G_{5,5})$, and $u_{\{1,2\}}v_1 \in M$, then

\begin{enumerate}
\item[(a)] there exists at least another saturated vertex in $V_1$, (either $u_4v_1, u_5v_1$ or both.
\item[(b)] suppose $G_a=G_{5,5} \backslash V_1$, and $M_a$ is some induced matching in $G_a$, then $|M_a| \neq \textrm{Max} (G_a)$.
\item[(c)] there exists some $e \in M$ with $e=u_iu_j$, such that $u_i,u_j \in V_5$.
\end{enumerate}

\end{lem}

\begin{proof} \begin{enumerate} \item[(a)] By Theorem \ref{thm2.2}, $|M|=6$. Suppose that $u_{\{1,2\}}v_1 \in M$ and that no other vertex in $V_1$ is saturated. Clearly, $u_1v_2$ and $u_2v_2$ can not be saturated. Therefore, there can only be $2$ saturated vertices in $V_2$. Suppose $u_{[4,5]}v_2 \in M$. Now we show that in $G_{5,5} \backslash \left\{V_1,V_2\right\}$, only $3$ edges belong to $M$. Let vertices $u_1v_3,u_1v_4, u_1v_5$ and $u_2v_3,u_2v_4,u_2v_5$, induce $G_b=G_{2,3} \subset G_{5,5}$. Now, $G_b$ has a maximum of four saturable vertices. Also let $G_c$ be a subgraph of $G_{5,5}$, induced by $u_3v_3,u_3v_4,u_3v_5$; $u_4v_4,u_4v_5$ and $u_5v_4,u_5v_5$. The subgraph $G_c$ also has a maximum of four saturable vertices. However, if we consider the positions $G_b$ and $G_c$, it is clear that at least a saturable vertex in $G_b$ is adjacent to a saturable vertex in $G_c$, which implies that $|V_s(G_b \cup G_c)| \leq 6$. Thus $|V_s(G_{5,5})| \leq 10$, and therefore $|M| \neq \textrm{Max} (G_{5,5})$ and hence a contradiction. By following similar argument, it can be seen also that $|M|=5$ if we consider $u_{\{3,4\}v_2} \in M$.
\item[(b)] Suppose that $G_a=G_{5,5} \backslash \{V_1\} \subset G_{5,5}$. Assume that $V_1$ contains only $3$ saturated vertices. There exists a saturated vertex $u_2 \in V_2$ such that given some vertex $u_1 \in V_1$, $v_1v_2 \in M$. Now, it is clear that $v_1v_2 \notin E(G_a) $ therefore, $|V_{sb}(G_a) \backslash v_2|=9$, implying that $|M_a|\leq 4$, while $\textrm{Max} (G_a)=5$.
\item[(c)] By (a) and (b) above, $|M_a|=4$. Suppose that $u_3v_{\{2,3\}} \notin M$, then, for $G_{5,5} \backslash \{V_1,V_2\}$, there must be at least $2$ saturated vertices on $V_5$. Let $G''=G_{5,2}=G_{5,5} \backslash \{V_1,V_2,V_5\}$ $\subset G_{5,5}$. Clearly $G_{5,2}$ will contain at most $6$ saturated vertices. Suppose that $u_a,u_b$ are saturated in $V_5$, and $u_au_b \notin E(G_{5,5})$, then $u_a$ and $u_b$ form $2$ edges with adjacent vertices $v_a,v_b \in V_4$. However, it can be seen that if this is so, there would be, at most, only $4$ saturable vertices in $G''$, apart from $v_a$ and $v_b$, and at least $1$ of which is an isolated vertex. Thus, there could only be $4$ saturated vertices in $G''$, which is a contradiction. Suppose $u_3v_{\{2,3\}} \in M$. It is easy to see by observation that it is impossible to have the matching in subgraph $G_f \subset G_{5,5}$, induced by $\left\{u_1v_3,u_1v_4,u_1v_5,u_2v_4,u_2v_5,u_3v_5,u_4v_4,u_4v_5,u_5v_3,u_5v_4,u_5v_5\right\} \in V(G_{5,5})$ since it contains at most $3$ edges in $M$ without any of them being made up of adjacent vertices in $V_5$.
    \end{enumerate}
\end{proof}
\begin{lem}\label{lem3.9a} Let $M$ be a maximum induced matching of $G_{5,9}$. If $v_{\{1,2\}}u_1 \in M,$ and there is at most one more saturated vertex $v_1$ on $V_1$. Then, $v_1=u_{5}v_1$.
\end{lem}
\begin{proof}
It is obvious that if $u_{\{1,2\}}v_1 \in M$, then $u_3v_1 \notin V_s(G_{5,9})$. Suppose that $v_1=u_4v_1$, then since it is the only lone saturated vertex on $V_1$, then, $u_4v_{\{1,2\}} \in M$. For $G_a$, some $G_{5,2}$ grid, induced by $V_1,V_2$, clearly, there is no other saturable vertex in $V_2$. Now, let $G_b=G_{5,9} \backslash G_a$. The subgrid $G_b$ is a $G_{5,7}$ grid and $\textrm{Max}(G_b)=8$. Thus, $|M|=10$, and hence, not $\textrm{Max}(G_{5,9})$. Hence a contradiction.
\end{proof}

\begin{lem}\label{lem3.9} Suppose that there are at most two saturated vertices on $V_1 \subset V(G_{5,5})$, and that they are adjacent. Then, there are at least $3$ saturated vertices on $V_5 \subset V(G_{5,5})$, two of which are adjacent.
\end{lem}
\begin{proof} Suppose $M$ is the maximal induced matching of $G_{5,5}$, and  that $u_au_b \in M$, where $u_jv_1,u_{j+1}v_1 \in U_1$. (It should be noted that since there are two saturated vertices in $V_1$ and are adjacent, then by Lemma \ref{lem3.8}(a), none of $j$ and $j+1$ is either $1$ or $5$). Suppose that $G_b \subset G_{5,5}$, induced by $V_2,V_3,V_4$ and that $V_5$ only contains $2$ saturated vertices $u_iv_5,u_{i+1}v_5$ and are adjacent, (meaning also that neither $i$ nor ${i+1}$ is $0$ or $5$.) Let $j \neq i$ and obviously, without loss of generality, set $i=2$, while $j$ clearly becomes $3$. Let $G_d, G_e \subset G_b$ be two subgraphs of $G_b$ induced by vertex sets $\left\{u_1v_2,u_1v_3,u_1v_4,u_2v_3,u_2v_4\right\}$ and $\left\{u_3v_3,u_4v_2,u_4v_3,u_5v_2,u_5v_3,u_5v_4\right\}$. $E(G_d), E(G_e)$ can only have a member each in $M$ and if $u_{\{3,4\}}v_3 \in M$, then no member of $E(G_d)$ belongs $M$. Then $|M|\leq 4$, which is a contradiction. If $i=j=1$, let $G_d$ be induced by $\left\{u_1v_2,u_1v_3,u_1v_4, u_2v_3\right\}$ and $G_e$ by $\left\{u_2v_3, u_4v_2, u_4v_3, u_4v_4,u_5v_2,u_5v_3,u_5v_4\right\}$. Following similar argument as above, $E(G_d)$ and $E(G_e)$ have maximum of $3$ members in $M$ and thus, $M$ is $5$, which is a contradiction. Suppose $V_5$ contains $3$ saturated vertices, such that none is adjacent to another. Clearly the saturated vertices are $u_1v_5,u_3v_5$ and $u_5v_5$. Since they are saturated, then $u_1v_{\{3,4\}},u_3v_{\{3_4\}},u_5v_{\{3,4\}} \in M$. Therefore by this, only two vertices on $V_3$ is saturable. Obviously $|M| \leq 5$.
\end{proof}

{\tiny{
\begin{center}
\pgfdeclarelayer{nodelayer}
\pgfdeclarelayer{edgelayer}
\pgfsetlayers{nodelayer,edgelayer}
\begin{tikzpicture}
	\begin{pgfonlayer}{nodelayer}
	
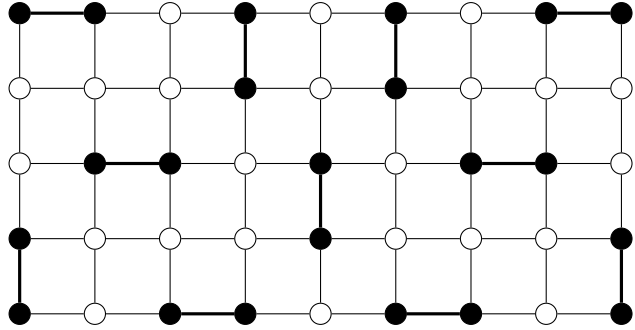
	\node [minimum size=0cm,]  at (-9,6.5) { A $G_{5,9}$ Grid with $Max(G_{5,9})=11$};

		\node [minimum size=0cm,draw,fill=black!,circle] (1) at (-13,7) {};
		\node [minimum size=0cm,draw,circle] (2) at (-12,7) {};
		\node [minimum size=0cm,draw,fill=black!,circle] (3) at (-11,7) {};
		\node [minimum size=0cm,draw,fill=black!,circle] (4) at (-10,7) {};
		\node [minimum size=0cm,draw,circle] (5) at (-9,7) {};
		\node [minimum size=0cm,draw,fill=black!,circle] (6) at (-8,7) {};
		\node [minimum size=0cm,draw,fill=black!,circle] (7) at (-7,7) {};
		\node [minimum size=0cm,draw,circle] (8) at (-6,7) {};
		\node [minimum size=0cm,draw,fill=black!,circle] (9) at (-5,7) {};
		\node [minimum size=0cm,draw,fill=black!,circle] (11) at (-13,8) {};
		\node [minimum size=0cm,draw,circle] (12) at (-12,8) {};
		\node [minimum size=0cm,draw,circle] (13) at (-11,8) {};
		\node [minimum size=0cm,draw,circle] (14) at (-10,8) {};
		\node [minimum size=0cm,draw,fill=black!,circle] (15) at (-9,8) {};
		\node [minimum size=0cm,draw,circle] (16) at (-8,8) {};
		\node [minimum size=0cm,draw,circle] (17) at (-7,8) {};
		\node [minimum size=0cm,draw,circle] (18) at (-6,8) {};
		\node [minimum size=0cm,draw,fill=black!,circle] (19) at (-5,8) {};
	  \node [minimum size=0cm,draw,circle] (21) at (-13,9) {};
		\node [minimum size=0cm,draw,fill=black!,circle] (22) at (-12,9) {};
		\node [minimum size=0cm,draw,fill=black!,circle] (23) at (-11,9) {};
		\node [minimum size=0cm,draw,circle] (24) at (-10,9) {};
		\node [minimum size=0cm,draw,fill=black!,circle] (25) at (-9,9) {};
		\node [minimum size=0cm,draw,circle] (26) at (-8,9) {};
		\node [minimum size=0cm,draw,fill=black!,circle] (27) at (-7,9) {};
		\node [minimum size=0cm,draw,fill=black!,circle] (28) at (-6,9) {};
		\node [minimum size=0cm,draw,circle] (29) at (-5,9) {};
		\node [minimum size=0cm,draw,circle] (31) at (-13,10) {};
		\node [minimum size=0cm,draw,circle] (32) at (-12,10) {};
		\node [minimum size=0cm,draw,circle] (33) at (-11,10) {};
		\node [minimum size=0cm,draw,fill=black!,circle] (34) at (-10,10) {};
		\node [minimum size=0cm,draw,circle] (35) at (-9,10) {};
		\node [minimum size=0cm,draw,fill=black!,circle] (36) at (-8,10) {};
		\node [minimum size=0cm,draw,circle] (37) at (-7,10) {};
		\node [minimum size=0cm,draw,circle] (38) at (-6,10) {};
		\node [minimum size=0cm,draw,circle] (39) at (-5,10) {};
		\node[minimum size=0cm,draw,fill=black!,circle] (41) at (-13,11) {};
		\node [minimum size=0cm,draw,fill=black!,circle] (42) at (-12,11) {};
		\node [minimum size=0cm,draw,circle] (43) at (-11,11) {};
		\node [minimum size=0cm,draw,fill=black!,circle] (44) at (-10,11) {};
		\node [minimum size=0cm,draw,circle] (45) at (-9,11) {};
		\node [minimum size=0cm,draw,fill=black!,circle] (46) at (-8,11) {};
		\node [minimum size=0cm,draw,circle] (47) at (-7,11) {};
		\node [minimum size=0cm,draw,fill=black!,circle] (48) at (-6,11) {};
		\node [minimum size=0cm,draw,fill=black!,circle] (49) at (-5,11) {};

		\end{pgfonlayer}
	\begin{pgfonlayer}{edgelayer}
		\draw [thin=1.00] (1) to (2);
		\draw [thin=1.00] (2) to (3);
		\draw [very thick=1.00] (3) to (4);
		\draw [thin=1.00] (4) to (5);
		\draw [thin=1.00] (5) to (6);
		\draw [very thick=1.00] (6) to (7);
		\draw [thin=1.00] (7) to (8);
		\draw [thin=1.00] (8) to (9);
		\draw [thin=1.00] (11) to (12);
		\draw [thin=1.00] (12) to (13);
		\draw [thin=1.00] (13) to (14);
		\draw [thin=1.00] (14) to (15);
		\draw [thin=1.00] (15) to (16);
		\draw [thin=1.00] (16) to (17);
		\draw [thin=1.00] (17) to (18);
		\draw [thin=1.00] (18) to (19);
		\draw [thin=1.00] (21) to (22);
		\draw [very thick=1.00] (22) to (23);
		\draw [thin=1.00] (23) to (24);
		\draw [thin=1.00] (24) to (25);
		\draw [thin=1.00] (25) to (26);
		\draw [thin=1.00] (26) to (27);
		\draw [very thick=1.00] (27) to (28);
		\draw [thin=1.00] (28) to (29);
		
		\draw [thin=1.00] (31) to (32);
		\draw [thin=1.00] (32) to (33);
		\draw [thin=1.00] (33) to (34);
		\draw [thin=1.00] (34) to (35);
		\draw [thin=1.00] (35) to (36);
		\draw [thin=1.00] (36) to (37);
		\draw [thin=1.00] (37) to (38);
		\draw [thin=1.00] (38) to (39);
		
		\draw [very thick=1.00] (41) to (42);
		\draw [thin=1.00] (42) to (43);
		\draw [thin=1.00] (43) to (44);
		\draw [thin=1.00] (44) to (45);
		\draw [thin=1.00] (45) to (46);
		\draw [thin=1.00] (46) to (47);
		\draw [thin=1.00] (47) to (48);
		\draw [very thick=1.00] (48) to (49);
		
		\draw [very thick=1.00] (1) to (11);
		\draw [thin=1.00] (11) to (21);
		\draw [thin=1.00] (21) to (31);
		\draw [thin=1.00] (31) to (41);
		
		\draw [thin=1.00] (2) to (12);
		\draw [thin=1.00] (12) to (22);
		\draw [thin=1.00] (22) to (32);
		\draw [thin=1.00] (32) to (42);
		
		\draw [thin=1.00] (3) to (13);
		\draw [thin=1.00] (13) to (23);
		\draw [thin=1.00] (23) to (33);
		\draw [thin=1.00] (33) to (43);
		
		\draw [thin=1.00] (4) to (14);
		\draw [thin=1.00] (14) to (24);
		\draw [thin=1.00] (24) to (34);
		\draw [very thick=1.00] (34) to (44);
		
		\draw [thin=1.00] (5) to (15);
		\draw [very thick=1.00] (15) to (25);
		\draw [thin=1.00] (25) to (35);
		\draw [thin=1.00] (35) to (45);
		
		\draw [thin=1.00] (6) to  (16);
		\draw [thin=1.00] (16) to (26);
		\draw [thin=1.00] (26) to (36);
		\draw [very thick=1.00] (36) to (46);
		
		\draw [thin=1.00] (7) to  (17);
		\draw [thin=1.00] (17) to (27);
		\draw [thin=1.00] (27) to (37);
		\draw [thin=1.00] (37) to (47);
		
		\draw [thin=1.00] (8) to  (18);
		\draw [thin=1.00] (18) to (28);
		\draw [thin=1.00] (28) to (38);
		\draw [thin=1.00] (38) to (48);
		
		\draw [very thick=1.00] (9) to (19);
		\draw [thin=1.00] (19) to (29);
		\draw [thin=1.00] (29) to (39);
		\draw [thin=1.00] (39) to (49);

	\end{pgfonlayer}
\end{tikzpicture}
\end{center}
}}

\begin{rem}\label{rem3.5} Let $G_a \subseteq G_{n,m}$ be a $G_{5,9}$ grid induced by $V_1, V_2, \cdots, V_9$, an induced matching $M$. From Lemma \ref{lem3.8} and Lemma \ref{lem3.9}, we see that if $u_{\{1,2\}}v_i \in M$, it is possible to have some $M_a \subseteq M$, for which $|M_a|=\textrm{Max} {(G_a)}$ as seen in the grid above. Also, from Lemma \ref{lem3.9}, since $V_1$ has $3$ saturated vertices with $2$ of them being adjacent, then for $M_a$ to be maximal,  $V_9$ will also have $3$ saturated vertices. It is easy to determine, however, that if this scheme extends to $m \geq 13$, then $|M| \neq \textrm{Max}(G_{n,m})$ since $|M| \leq 11+4j+5k$, $j \geq k$, and $j-k \leq 1$.
\end{rem}

\begin{lem}\label{lem3.10} Let $M$ be a matching of $G_{5,m}$, $m\equiv 1 \mod 4$, $m \geq 13$, with $u_{\{1,2\}}v_1 \in M$, then $|M| \neq \textrm{Max} (G_{5,m})$.
\end{lem}

\begin{proof}
 Let $m-10=q$. Clearly, $q \equiv 3 \mod 4$. For $q \geq 7$, let $G_{5,q}= G_a$, induced by $V_{11},V_{12}, \cdots V_m$. We already know that $|V_{sb}(G_a)|=\frac{5q+1}{2}$, while $|V_s(G_a)|=2\left(\frac{5(q-1)-2}{4}\right)$. Therefore, there are $2$ saturable vertices, say, $v_1,v_2$ on $V(G_a)$. Suppose that $v_1,v_2 \in V_{11}$. Note that $v_1,v_2$ are not adjacent, else $v_1v_2 \in M$. Let $G_b \subset G_{5,m}$ be a $G_{5,9}$, induced by $V_1,V_2,\cdots, V_9$. Suppose that $u_{\{1,2\}}v_1 \in M $, then there are two adjacent saturated vertices in $V_9$ and at least a saturated boundary vertex,  by Lemma \ref{lem3.9} and Remark \ref{rem3.5}. 
  This implies that there are two adjacent saturable vertices on $V_{10}$, say $u_1,u_2$. Now, there could be at most one edge in $M$ from $\left\{u_1,u_2,v_1,v_2\right\}$. Hence, $|M|=\textrm{Max} (G_b)+\textrm{Max} (G_a)+1=\frac{45+5q}{4}$. Since $q=m-10$, we have $|M|= \frac{5(m-1)}{4}$, which is less than $\textrm{Max} (G_{5,m})$ by an edge, and hence a contradiction. For $q=3$, the result is similar to that obtained by careful observation of the positions $u_1,u_2$ and the possible isolated vertex or vertices on $G_a$.
\end{proof}

\begin{corollary}\label{cor3.11} Let $M$ be a matching of $G_{5,m}$, $m \geq 13$, with $u_{\{1,2\}}v_i \in M$, $1 \leq i \leq m$, then $|M| \neq \textrm{Max} (G_{5,m})$.
\end{corollary}

\begin{lem}\label{lem3.11} Let $U_1 \subset V(G_{5,m})$. Then there are at least $\frac{m+1}{2}$ saturable vertices in $U_1$.
\end{lem}

\begin{proof} This follows similar arguments as in the proof of Theorem \ref{thm3.4}.

\end{proof}
\begin{rem}\label{rem3.6} We note that for $m=4k+1$, $\frac{m+1}{2}$ is odd, also from Lemma \ref{lem3.10}, the number of saturated vertices on $U_1$ is even and there cannot be any isolated vertex on $U_1$ if $M$ is a maximal induced matching of $G_{5,m}$. Therefore, for $G_1 \subset G_{5,m}$, induced by $U_1$, there exists at least $k+1$ edges of $G_1$ in $M$, for $m \geq 13$. For $m=9$, there are at least $k$ edges in $G_1$ as seen in the last figure.
\end{rem}
\begin{thm}\label{thm24}
Let $G_{n,m}$ be a grid with $m \equiv 1 \mod 4$, $m\geq 13$ and let $M$ be the maximum induced matching of $G_{n,m}$.  Then
$$ \textrm{Max}(G_{n,m}) \leq \left \{ \begin{array}{ll} \left\lfloor \frac{2mn-m-1}{8}\right\rfloor & \  \textrm{if} \ n \equiv 1 \mod 4; \\ \ \\
\left\lfloor \frac{2mn-m+3}{8}\right\rfloor & \ \textrm{if} \ n \equiv 3 \mod 4
\end{array} \right. $$



\end{thm}

\begin{proof} Let $n=4l+1$, $l$ a positive integer, and let $r=n-5$, that is, $r \equiv 0 \mod 4$. Suppose $G_a$ is a $G_{r,m}$ induced by $V_1,V_2, \cdots, V_r$. By Theorem \ref{lem3.11}, at least there are $\frac{m+1}{2}$, saturated vertices on $U_r$ and these are $u_rv_2,u_rv_4, \cdots, u_rv_{m-1}$. Let $G_b$ be a $G_{5, m}$ grid, induced by $V_{r+1}, V_{r+2}, \cdots, V_n$ and $G_c$ be a $G_{1,m}$ grid, induced by $U_{r+1}$. By Remark \ref{rem3.6}, there are $k+1$ edges of $E(G_c)$ in $M$. Clearly a saturated vertex induced by $M$ on each of the $k+1$ edges is adjacent to some saturated vertices in $U_r$, implying that only one of the two saturated vertices belongs to $V_{sb}(G_{n, m})$. Hence, $V_{sb}(G_{n, m}) \leq \frac{nm+1}{2}-k+1$. Now, $k+1 = \frac{m-1}{4}+1$. Therefore, $V_{sb}(G_{n,m}) \leq \frac{2nm-m-1}{4}$ and hence, $|M| \leq \left\lfloor \frac{2nm-m-1}{8}\right\rfloor$. For $n=4l+3$, set $s=n-3$, that is, $s \equiv 0 \mod 4$. By Remark \ref{rem3.6}, and following the argument above, $V_{sb}(G_{n,m})=\frac{mn+1}{2}-k$ with $k=\frac{m-1}{4}$, $\textrm{Max}(G_{n,m}) \leq \left\lfloor \frac{2mn-m+3}{8}\right\rfloor$.

\end{proof}
\begin{rem} For the grid $G_{n,m}$, it should be noted that the results in the last Theorem depends mainly on the value of $m$.
\end{rem}

\begin{rem}Following similar argument as Theorem \ref{thm24}, for $m=9$, $\textrm{Max}(G_{n,9}) \leq  \lfloor \frac{18n-6}{8}\rfloor.$
\end{rem}

\section*{Acknowledgements} We would like to thank an anonymous referee for very helpful comments and suggestions.

%
%
%
%






\begin{thebibliography}{99}
		
		\bibitem{C1} K. Cameron. Induced matching in intersection graphs. \emph{Discrete Math.} \textbf{278} (2004), 1--9.
		
		\bibitem{CST1} K. Cameron, R. Sritharan and Y. Tang. Finding a maximum induced matching in weakly chordal graphs. \emph{Discrete Math.} \textbf{266} (2003), 133-142.
		
		 \bibitem{DD1} K. K. Dabrowski, M. Demange and V. Lozin. New results on maximum induced matchings in bipartite graphs and beyond. \emph{Theoretical Computer Science,} \textbf{478} (2013), 33--40.
		
		\bibitem{E1} J. Edmonds, Paths, trees, and flowers, \emph{Canad. J. Math.} \textbf{17} (1965), 449--467 .
		
		  \bibitem{GL1} M. C. Golumbic and R.C Laskar, Irredundancy in circular arc graphs. \emph{Discrete Applied Math.} \textbf{44} (2013), 79--89.
		
		  \bibitem{J1} F. Joos, Induced matching in graphs of bounded maximum degree. \emph{arXiv:} (2014) 1406.2440.
		
			\bibitem{RMG1} R. Marinescu-Ghemeci,  Maximum induced matchings in grids. \emph{Optimization Theory, Decision Making and Operations Research Applications. Springer New York}, (2013)  177--187.
		
		  \bibitem{N1} V. H. Nguyen,  Induced matching in graphs of degree at most 4. \emph{arXiv:} (2014) 1407.7604.
		
		  \bibitem{SV1} L. J. Stockmeyer and V.V. Vazirani. $NP-$completeness of some generalizations of the maximum matching problem. \emph{Information Processing Letters}, \textbf{15}(1) (1982), 14--19.
		
		  \bibitem{Z1} M. Zito,   \emph{Induced matching in regular graphs and trees}. Lecture Notes in Computer Sci. 1665, Springer, Berlin, 1999
		
	
	
	
	
	
	
	
	
	
	
	
	
	
		
		
	\end{thebibliography}
\end{document}